\documentclass[]{amsart}
\title[Residual properties of mapping tori]{Polynomial maps over $p$-adics and
residual properties of mapping tori of group endomorphisms}

\author{Alexander Borisov}
\address{Department of Mathematics, University of Pittsburgh}
\email{borisov@pitt.edu}
\author{Mark Sapir}
\address{Department of Mathematics, Vanderbilt University}
\email{mark.sapir@vanderbilt.edu}
\thanks{The research of the
first author was supported in part by the NSA grants H98230-06-1-0034 and H98230-08-1-0129. The work of the second author
was  supported in part by the NSF grant DMS 0700811 and by a BSF
(USA-Israeli) grant.}

\long\def\comment#1\endcomment{}

\newcommand{\la}{\langle}
\newcommand{\ra}{\rangle}

\newcommand{\alg}{{\mathrm{alg}}}

\newcommand{\iv}{^{-1}}

\newcommand{\F}{\mathbb{F}}
\newcommand{\Fr}{\hbox{Fr}}
\newcommand{\HNN}{\hbox{HNN}}
\newcommand{\Q}{{\mathbb Q}}

\newcommand{\Z}{{\mathbb Z}}
\newcommand{\N}{{\mathbb N}}

\newcommand{\A}{{\rm A}}
\newcommand{\OO}{{\mathcal O}}
\newcommand{\SL}{\hbox{SL}}

\newcommand{\PGL}{\hbox{PGL}}
\newcommand{\Spec}{{\rm Spec}}

\newcommand{\adj}{{\hbox{adj}}}
\begin{document}

\begin{abstract} We continue our study of residual properties of mapping tori of free group endomorphisms. In this paper, we prove that each of  these groups are virtually residually (finite $p$)-groups for all but finitely many primes$p$. The method involves further studies of polynomial maps over finite fields and $p$-adic completions of number fields.
\end{abstract}

\maketitle

\renewcommand{\theequation}{\thesection.\arabic{equation}}
\newtheorem{theo}{\quad Theorem}[section]
\newtheorem{lemma}[theo]{\quad Lemma}
\newtheorem{cy}[theo]{\quad Corollary}
\newtheorem{df}[theo]{\quad Definition}
\newtheorem{rk}[theo]{\quad Remark}
\newtheorem{prop}[theo]{\quad Proposition}
\newtheorem{prob}[theo]{\quad Problem}
\newtheorem{conj}[theo]{\quad Conjecture}

\section{Introduction}

This paper is a continuation of paper \cite{BS} where we proved that
for every linear finitely generated group $G$ and any injective
endomorphism $\phi$ of $G$, the mapping torus of $\phi$ is
residually finite.  The mapping torus of $\phi$ is the following
{\em ascending} HNN extension of $G$:

$\HNN_\phi(G)=\la G, t\mid txt\iv =\phi(x)\ra$ where $x$ runs over a (finite) generating set of $G$.

Probably the most important mapping tori are mapping tori of
endomorphisms of free groups $F_k$. These groups appear frequently
as fundamental groups of hyperbolic 3-manifolds (in fact there is a
conjecture that all fundamental groups of hyperbolic 3-manifolds are
virtually mapping tori of free group automorphisms). Also it is
proved in \cite{KS}, that with probability tending to 1 as $n\to
\infty$, every 1-related group with three or more generators and
relator of length $n$ is embeddable into the mapping torus of a free
group endomorphism (and so it is residually finite by \cite[Theorem
1.6]{BS} and coherent by \cite{FH}). For 1-related groups with 2
generators, it is not known whether they are almost surely inside
mapping tori of free group endomorphisms, but the computer
experiments from \cite{BS} and \cite{DT} show that the probability
of that should be at least $.94$.

A fairly comprehensive survey of our knowledge about mapping tori of
free group endomorphisms before \cite{DS,BS} is in \cite{Kap}. The
paper \cite{Kap} ends with several open problems. The first problem
asks whether $\HNN_\phi(F_k)$ is residually finite or linear.
\cite{BS} answers positively the first part of the question and
\cite{DS} answers negatively the second part: it turned out that the
group $\langle a,b,t\mid tat^{-1} =a^2, tbt^{-1} =b^2\rangle$ is not
linear. After \cite{BS} and \cite{DS}, a natural question arises: do
groups $\HNN_\phi(F_k)$ possess stronger residual properties than
simple residual finiteness? For example, are they residually finite
nilpotent groups, etc? The answer to that question is ``no":
consider any endomorphism $\phi$ of $F_k$ that maps $F_k$ into the
derived subgroup $[F_k,F_k]$. Then every solvable homomorphic image
of $H=\HNN_\phi(F_k)$ is cyclic.

Indeed, let $H=\langle x_1,...,x_k,t\mid tx_it^{-1}=\phi(x_i),
i=1,...,k\rangle$. Then $\phi(x_i)$ is a product of commutators in
$F_k$ for every $i=1,...,k$. Let $H'$ be a non-cyclic solvable
homomorphic image of $H$. Let $\bar x_i$ be the image of $x_i$ and
$\bar t$ be the image of $t$ in $H'$. Then for every $i$, the element
$\bar t \bar x_i\bar t^{-1}$ is a product of commutators in the
subgroup $X$ of $H'$ generated by $t^n\bar x_1t^{-n}, ..., t^n\bar
x_kt^{-n}$, $n\in \Z$. Note that $H'$ is an extension of $X$ by a
cyclic group $\langle \bar t\rangle$. Hence $\bar t^n\bar x_i\bar t^{-n}$ is in the
derived subgroup of $X$ for every $n,i$. Hence $X$ coincides with
$[X,X]$. Since $X$ is solvable, $X$ must be trivial, and $H'$ must
be cyclic, a contradiction.

Note that many linear groups (for example,  $\SL_n(\Z)$ for $n>2$) are also not
residually solvable (by Margulis's normal subgroup theorem), but all
finitely generated linear groups  are {\em virtually} residually
solvable and even virtually residually  (finite $p$)-groups for all
but finitely many primes $p$ \cite{Mal}. That means they have finite
index subgroups that are residually (finite $p$)-groups.  In this
paper, we shall prove that all mapping tori of finitely generated
free group endomorphisms enjoy this property too.

\begin{theo}\label{1} Every mapping torus of a finitely generated free group endomorphism is virtually residually (finite $p$)-group for every sufficiently large $p$.
\end{theo}

We shall illustrate the ideas of the proof by the following example (which was the motivating example for our work).

Let $H=\langle a,b,t\mid tat^{-1}=ab, tbt^{-1}=ba\rangle$. It is not
difficult to prove (almost as above) that this group is not
residually solvable. We shall prove that it has a subgroup of finite
index which is residually (finite $5$)-group. Consider two matrixes
$A=\left(\begin{array}{cc}5 & 2\\ 2 & 1\end{array}\right)$,
$B=\left(\begin{array}{cc}1 & 2\\ 2 & 5\end{array}\right)$. These
two matrices generate a free subgroup of $\SL_2(\Z)$ (which can be
easily proved either by noticing that it is a subgroup of Sanov's
free subgroup of $PSL_2(\Z)$, or by just looking at the subgroup
generated by the corresponding Frobenius transformations and using
the ping-pong argument). For every group $G$ consider the
map $\phi_G\colon  G\times G\to G\times G$ given by $\phi(U,V)=(UV,
VU)$. It is easy to check that if $G_1=\SL_2(\Z/5\Z)$, then the pair
$(A \mod 5, B \mod 5)$ is a periodic point for $\phi_{G_1}$ with
period 6. Replacing $5$ by $5^2=25$, and considering the group
$G_2=\SL_2(\Z/25\Z)$,  one can compute that $(A \mod 25, B \mod 25)$
is a periodic point of $\phi_{G_2}$ with period $5*6=30$. By
induction or by using Theorem  \ref{MainAlg} below, one can prove
that for $G_n=\SL_2(Z/5^n\Z)$, the pair $(A \mod 5^n, B \mod 5^n)$
is a periodic point for $\phi_{G_n}$ with period $q_n=6*5^{n-1}$. It
is observed in \cite[Lemma 2.2]{BS}, that $\HNN_\phi(F_k)$ has
homomorphisms $\nu_n$ into the wreath product $G_n\wr \Z/q_n\Z$
which is the semi-direct product of $G_n^{q_n}$ and the cyclic group
$\Z/q_n\Z$ where the cyclic group acts on the direct power by the
cyclic shift. The homomorphism $\nu_m$ takes $t$ to the generator
$t_n$ of $\Z/q_n\Z$, $a$ maps to $a_n=(A \mod 5^n, \phi_{G_n}(A \mod
5^n), ...,  \phi_{G_n}^{q_n-1}(A \mod 5^n))$, $b$ maps to $b_n=(B
\mod 5^n, \phi_{G_n}(B \mod 5^n), ..., \phi_{G_n}^{q_n-1}(B \mod
5^n))$. There exists a natural homomorphism  $\mu_n\colon G_n\to
G_1$. The image $\langle a_n, b_n\rangle$ of $F_2=\langle
a,b\rangle$ under $\nu_n\mu_n$ is inside the direct power
$G_1^{q_n}$, hence $\nu_n\mu_n(F_2)$ is a 2-generated group in the
variety\footnote{The creators of the theory of varieties of groups
(G. Birkhoff mostly) did not anticipate in the 1930s the situation when the word
``variety" would be used in the same paper in two different senses:
as an algebraic variety and as a variety of algebraic systems, i.e.
a class of algebraic systems closed under taking cartesian products,
subsystems and homomorphic images. This is one of the very few
papers (if not the only paper) where the term is used in both
senses.} of groups generated by the finite group $G_1$. It is well
known that then $\nu_n\mu_n(F_2)$ has order bounded by some constant
$M$. Then the centralizer of $\mu_n\nu_n(F_k)$ in $\mu_n\nu_n(H)$
has index at most $M_1=M^M$ in $\mu_n\nu_n(H)$. Since $\Z/q_n\Z$ has
a 5-subgroup of index 6, $\mu_n\nu_n(H)$ has a 5-subgroup of index
at most some constant $M_2$ (independent of $n$). Since by
\cite[Lemma 2.2]{BS}, all elements of the group $H$ are separated by
homomorphisms $\nu_n$, $H$ has a subgroup of index at most $M_2$
which is residually (finite $5$)-group.

One can follow the proof of Theorem \ref{1} below (see Section \ref{sect3}) in order to show that $H$ is virtually residually (finite $p$)-group for almost all $p\ne 5$. The differences with the proof above are the following
\begin{itemize}

\item a pair of matrices $(A,B)$ such that $(A \mod p, B \mod p)$ is
periodic for $\phi_{G_1}$ is found not in $G_1=\SL_2(\Z/5\Z)$, but
in $G_1=\SL_2(\OO/p\OO)$ where $\OO$ is the ring of integers of some
finite extension of $\Q$ unramified at $p$ (we use the main
algebro-geometric result \cite[Theorem 3.2]{BS} for this);
\item in order to find such a pair of matrices with the  additional property that  $(A, B)$ generate a free subgroup, we use a strong result of Breuillard and  Gelander \cite{BG} about dense free subgroups of Lie groups;  the matrices $A, B$ are found not in $\SL_2(\OO)$ but in the $p$-adic completion of that group. Here we use a new result about polynomial maps over $p$-adics proved in Section \ref{sect2} below, see Theorem \ref{MainAlg}.
\end{itemize}

\section{Polynomial maps over $p$-adics}
\label{sect2}

For a prime power $q=p^k,$ denote by $\F_q$ the field of $q$
elements. Denote by $\F_q^{alg}$ its algebraic closure
($\F_q^{alg}=\F_p^{alg}$).

An affine algebraic variety over $\F_q^{alg}$ is a subset of
$(\F_q^{alg})^{n}$ consisting of all common roots of some ideal
$I\subset \F_q^{alg}[x_1,...,x_n].$ By Hilbert's Nullstellensatz, we
can assume that $I$ is radical: for every $f\in
\F_q^{alg}[x_1,...,x_n]$ and $N\in {\mathbb{N}},$ if $f^N\in I$ then
$f\in I.$

If $X$ is an affine algebraic variety over $\F_q^{alg},$ its field
of definition is the smallest subfield of $\F_q^{alg}$ containing
all coefficients of some set of generators of $I.$ An affine
algebraic variety over $\F_q$ is a variety with the field of
definition is $\F_q$ or its proper subfield. Equivalently, $X$ is
defined over $\F_{p^k}$ if and only if $\Fr ^k (X)=X,$ where $\Fr$
is the geometric Frobenius self-map of the affine space: $\Fr
(x_1,x_2,...,x_n)=(x_1^p, x_2^p,...,x_n^p).$

Note that in algebraic geometry one usually wants to understand the
structure of $X$ independent of its embedding into the affine space,
as the affine scheme associated to the ring $\F_q[x_1,...,x_n]/I.$
However all varieties in this paper naturally appear as subvarieties
of the fixed affine space. We will not be dealing with scheme points
of $X$, but rather with its geometric points, which are just the
points $(x_1,x_2,...,x_n)$ of the affine $n-$space over
$\F_q^{alg}$,  contained in $X$ (or, more precisely, its base change
to $\F_q^{alg}.$

Let, as usual, $\Z_p, \Q_p$ be the $p$-adic completion of the ring of
integers and rational numbers, respectively. For $q=p^k$ let $\Q _q$
be the unrafimied extension of $Q_p$ such that the residue field
$Z_q/pZ_q$  of its ring of integers $Z_q$ is isomorphic to $\F_q$ (see, e.g. \cite[page 143, example 4.18]{MP}).
Let ${\A}^n$ be the affine space over $\Z$ of dimension $n$, that is
${\rm{Spec}} {\mathbb Z}[x_1,...,x_n]$.

\begin{df} The following are standard terms in algebraic geometry, adapted for our purposes.

\begin{itemize}
\item An algebraic variety is called {\em geometrically irreducible} if
its base change to the variety over the algebraic closure of the
ground field is irreducible. This means that the ring
$\F_q^{alg}[x_1,...,x_n]/I$ is a domain.
\item For a variety $X$ as above, a polynomial self-map $\Phi$ of $X$ is a polynomial self-map of an ambient affine space that preserves $X$. In coordinates, $\Phi$ is given by polynomials:
$$\Phi (x_1,x_2, ... , x_n)=(f_1(x_1,...,x_n),f_2(x_1,...,x_n),...,f_n(x_1,...,x_n))$$
The field of definition of $\Phi$ is the subfield of $\F_q^{\alg}$ generated by the coefficients of $f_i$. $\Phi$ preserves $X$ when for every $g\in  \F_q^{alg}[x_1,...,x_n]$ that belongs to $I$, the formal composition $g\circ \Phi$ also belongs to $I$ (I is assumed to be radical). We will only consider the polynomial self-maps.

The self-map $\Phi$ is called {\em dominant} if its image (over $\F_q^{alg}$) is Zariski dense in $V$.

\item A dominant self-map $\Phi$ is called {\em separable} if the corresponding extension of the fields of rational function is separable.
\item Zariski tangent space at a point $(x_1,x_2,...,x_n)\in X,$ defined over $\F_Q,$ is the $\F_Q-$vector space, dual to the space $m/m^2,$ where $m$ is the ideal in $\F_Q[x_1.x_2,...x_n],$ consisting of polynomials that vanish on   $(x_1,x_2,...,x_n)$.
\item A self-map $\Phi\colon V(\F_q)\to V(\F_q)$ is called {\em unramified} at a point $x\in V$ if the map $\Phi_{*}$ on Zariski tangent space of $x$ is invertible. (Note that a separable dominant self-map is unramified at all $x$ in some Zariski open subset of $V$).
\item A geometric point $x$ of an algebraic variety $V$ is called {\em smooth} if the local ring is regular. Alternatively, the point is smooth if the Zariski tangent space at it has the same dimension as the variety.
\item The {\em degree} of $\Phi$ is the number of geometric points in the preimage of a  generic point of $V$.

\end{itemize}
\end{df}

The following statement is proved in \cite{BS}.

\begin{lemma}(\cite[page 349]{BS})\label{lm1}
Suppose that $\Phi\colon {\A}^n\to {\A}^n$ is a polynomial map
defined over any algebraically closed field. Denote by $V$ the
Zariski closure of $\Phi^n({\A}^n)$. Then $V$ is geometrically
irreducible and the map $\Phi|_V\colon V\to V$ is dominant.
\end{lemma}

Suppose now that we have a polynomial map $\Phi$ from ${\A}^n$ to itself, defined over $\Z$. Applying the above lemma to the map over ${\Q}^{\alg},$ we get some subvarity $V$ of ${\A}^n$. Even though it is a priori defined only over ${\Q}^{\alg},$ its field of definition is ${\Q},$ because it is fixed by the absolute Galois group of ${\Q}.$ Reducing $\Phi$ modulo $p,$ we get a polynomial self-map of ${\A}^n$ over $\F_p,$ and thus over $\F_p^{\alg},$ the algebraic closure of $\F_p$. This in turn produces a subvarity $V_p$ of ${\A}^n(\F_p^{\alg}).$ Naturally, we would like to relate $V_p$ to $V$ for large enough $p.$ In order to do this, we need to construct a model of $V$, that is a scheme over $\Spec {\Z},$ such that $V(\Q)$ is its generic fiber. Fortunately, a natural model exists in our situation. To describe it, let's consider $V$, and $V_p,$ from a commutative algebra perspective.

Suppose $\Phi=(\Phi_1, \Phi_2, ..., \Phi_n)$ is an ordered set of polynomials in ${\Z} [x_1,x_2,...,x_n] $, and the varieties $V$ and $V_p$ are defined as above. Abusing the notation a little bit, denote by $\Phi ^*$  the corresponding ring homomorphism, for $n-$variable polynomials over any ring. Then the prime ideal $I(V)$ of ${\Q}^{\alg} [x_1,x_2,...,x_n] ,$  defining $V, $ consists of polynomials $f \in {\Q}^{\alg} [x_1,x_2,...,x_n] $ such that $f(\Phi^n) =0$, where $\Phi ^n$ is the $n-th$ composition power of $\Phi .$ This is, in other words, the pullback of the prime ideal $\{0\}$ of ${\Q}^{\alg} [x_1,x_2,...,x_n]$ by the homomorphism $(\Phi^n)^*=(\Phi^*)^n.$ As a pullback of a prime ideal, it is also prime, in particular radical. Similarly, the ideal  $I(V_p)$ of $V_p$ in ${\F_p}^{\alg} [x_1,x_2,...,x_n] $ consists of polynomials $f \in {\F_p}^{\alg} [x_1,x_2,...,x_n] $ that  $f((\Phi_p)^n) =0$, where $\Phi _p$ is the reduction of $\Phi$ modulo $p.$

Now we define a model of $V$ as an affine scheme over $\Z$ which is a subscheme of ${\Spec {\Z}}^n$ defined by the ideal $I(V(\Z)) \in \Z [x_1,x_2,...,x_n],$ consisting of all $f$ that $f(\Phi ^n)=0$ in $I \in \Z [x_1,x_2,...,x_n].$

\begin{lemma} The ideal $I(V(\Z))$ is prime. It is finitely generated and its generators when considered over $\Q^{\alg}$ generate the ideal $I(V).$
\end{lemma}

\proof The ideal $I(V(\Z))$ is prime because it is a pullback of the prime ideal $\{0\},$ like in the argument above. It is finitely generated because the ring $\Z [x_1,x_2,...,x_n]$ is Noetherian. For the last statement, it is obvious that all elements of $I(V(\Z))$ belong to $I(V).$ Suppose $f$ is any element of $I(V)\subset {\Q}^{\alg} [x_1,...,x_n].$ Because all coefficients of $\Phi$ are integers, all conjugates of $f$ also belong to $I(V)\subset {\Q}^{\alg} [x_1,...,x_n].$  Consider a ${\Q}^{\alg}-$vector subspace of a space of $n-$variable polynomials of large enough degree, generated by all these conjugates of $f$. This subspace is invariant under the natural action of the absolute Galois group of $\Q,$ thus it has a basis consisting of polynomials with rational coefficients. For every element in this basis some non-zero integer multiple of it
belongs to $\Z ^{\alg} [x_1,x_2,...,x_n].$ It is obvious that it must belong to $I(V(\Z)),$ which implies that $f$ is a linear combination with algebraic coefficients of elements from $I(V(\Z))$.

We would like to say that for all $p$ the reductions of the generators of $I(V(\Z))$ modulo $p$ generate the ideal $I(V(\Z/p\Z)),$ but this is in general not true, as the following example shows.

{\bf Example.} Suppose $n=2$, denote the coordinates by $x$ and $y$. Suppose $\Phi (x,y)=(x,5y)$. Then $\Phi ^2 (x,y)= (x, 25y).$ In characteristic zero this is an invertible linear map, so the ideal $I(V(\Z))$ is zero. However, modulo $p=5$ this map is not invertible and $I(V(\Z/p\Z))$ is the principal ideal, generated by the polynomial $x$.

The following theorem is very important. It shows that for all sufficiently large primes $p$ the reductions of the generators of $I(V(\Z))$ modulo $p$ generate the ideal $I(V(\Z/p\Z)).$

\begin{theo} Suppose $\Phi \in ({\Z}[x_1,...x_n])^n$ is a polynomial automorphism in $n$ variable with integer coefficients. Define $I(V(\Z))$ and $I(V(\Z/p\Z))$ are as above. Then there exists some natural number $B,$ such that for all $p>B$ the reductions of the generators of $I(V(\Z))$ modulo $p$ do generate the ideal $I(V(\Z/p\Z)).$
\end{theo}

\proof
For every $p$ we have the inclusion of the ideals in $\F_p[x_1,x_2,...,x_n]$: $I(V(\F_p))$ contains the reduction modulo $p$ of the ideal $I(V(\Z))$, to be denoted by $I_p$. By a very general result on fibers of algebraic morphisms of schemes (cf. \cite{Eisenbud}, Thm. 14.8), all irreducible components of $\Spec (\F_p [x_1,...,x_n] / I_p)$ have the dimension equal to the dimension of $\Q[x_1,...x_n]/I(V(Q))$. Because $I(V(\F_p))$ is prime, it is enough to show that the dimension of $\Spec \F_p[x_1,...,x_n]/I(V(\F_p))$ is also the same for all sufficiently large $p$. The key idea of this argument is the following construction.

Suppose $\F$ is any field, and $M$ is a polynomial self-map of the affine space $\Spec \F[x_1,x_2,...x_n].$ The map $M$ is given in coordinates by its components $M_i,$ which are polynomials in  $\F[x_1,x_2,...x_n].$ Consider its formal Jacobian matrix.
$J=(\frac{ \partial M_i}{\partial x_j})_{i,j=1,...n}$. The entries are in $\F[x_1,x_2,...x_n],$ so the matrix is well-defined as an $n\times n$ matrix $T$ over the field $\F(x_1,...,x_n)$.
We define the ideal $I$ in $\F[x_1,...x_n]$ as a pullback of the zero ideal by $M$.

\begin{lemma} In the above notation, if the map $M$ is separable, then the dimension of $\Spec \F[x_1,...,x_n]/I$ equals the rank of the matrix $T$.
\end{lemma}
\proof This is (the algebraic version of) the classical inverse function theorem, see e.g. \cite{Eisenbud}, Chapter 16  for the reference. \endproof

We now apply the above lemma to $M=\Phi ^n$ for various fields $\F$. For  $\F=\Q$ the rank of $T$ is the dimension of $V(\Q)$, let's denote it by $r.$ The entries of $T$ over $\Q$ are actually polynomials with integer coefficients. For some $r\times r$ minor of $T$ the determinant is a non-zero polynomial with integer coefficients. For all sufficiently large primes $p$ its reduction modulo $p$ is also non-zero. So for these $p$ the rank of the Jacobian matrix for $\F=\F_p$ is at least $r.$ Thus the dimension of $\Spec \F_p[x_1,...,x_n]/I(V(\F_p))$ is at least $r$. On the other hand, any $(r+1)$ rows of $T$ over $\Q$ are linearly dependent over $\Q(x_1,...x_n).$ Multiplying by a suitable polynomial in $\Z[x_1,...x_n]$, we get non-trivial linear combination of rows of $T$ with coefficients in $\Z[x_1,...x_n]$ that equals to zero. Since $T$ for $\F=\F_p$ is just the reduction of $T$ for $\Z$ modulo $p$, for all sufficiently large $p$ all of these linear combinations are still non-trivial. So the rank of $T$ is exactly $r$ for all sufficiently large $p.$ For such $p$ we have $\Spec I_p$ a possibly reducible subvariety of an irreducible variety $V(\F_p),$ and every component of $\Spec I_p$ has the same dimension $r$ as $V(\F_p).$  This implies that the varieties are the same, so the ideals $I(V(p))$ and $I_p$ are equal, for any sufficiently large $p$.

\endproof

From now on, the prime $p$ will always be sufficiently large so that the conclusion of Theorem 2.4 holds. Let us denote by $\pi_p$ (or just $\pi,$ when $p$ is  the reduction map from $\Z$ to $\Z / p \Z$ and, abusing the notation a bit, all other reduction modulo $p$ maps. Recall, that for any $q=p^k$ we denoted by $\Z_q$ the unramified extension of the $p-$adic integers $\Z_p$ with the residue field $\F_q.$ With the above convention, we will denote by $\pi_p$ the reduction map from $\Z_q$ to $\F_q,$ as well as the reduction map from $\Z_q[x_1,...x_n]$ to $\F_q[x_1,...x_n].$

In what follows, unless otherwise specified, we will denote by upper case letters objects in characteristic zero, and by corresponding lower case letters the objects in characteristic $p.$ For example, if $T\in \Z_q[x_1,...x_n]$ then $t\in \F_q[x_1,...x_n]$ is its reduction modulo $p:$ $t=\pi_p (T).$

\begin{lemma} For $p$ as above, for any $q=p^k$ define $I(V(\Z_q))\subset \Z_q[x_1,...x_n]$ and $I(V(\F_q)) \subset \F_q[x_1,...x_n]$ as before. Then
$$I(V(\F_q)) = \pi_p (I (V(\Z_q)))$$
\end{lemma}

\proof
We need to prove two inclusions, one of which is easy. Indeed, if $g=\pi_p(G), $ $G\in I (V(\Z_q)),$ then $G(\Phi^n) =0,$ so $G(\Phi^n) =0 \mod p,$ so $g\in I(V(\F_q)).$

Now we would like to prove that if $g\in I(V(\F_q)),$ there exists $G\in I(V(\Z_q))$ such that $g=\pi_p(G).$ Choose a set of generators $\{G_1,G_2,...,G_m\}$ of the ideal $I(V(\Z)),$ as the ideal in $\Z[x_1,...,x_n].$ By assumption, $g_i=\pi_p (G_i)$ generate $I(V(\F_p))$ as the ideal in $\F_p[x_1,...x_n].$ By the same argument as in the Lemma 2.3, this implies that they generate $I(V(\F_q))$ in $\F_q[x_1,...x_n].$

So, modulo $p$, $\pi_p(G)=\sum \limits_{(\alpha,\beta)}c_{(\alpha,\beta)}  g_1^{\alpha_1}\cdot ...\cdot g_m^{\alpha_m} \cdot x_1^{\beta_1}\cdot ... \cdot x_n^{\beta_n},$ where $\alpha$ and $\beta$ are multi-indices. Lifting $c_{(\alpha,\beta)}$ to $C_{(\alpha,\beta)},$ we get
$$G=\sum \limits_{(\alpha,\beta)}C_{(\alpha,\beta)}  G_1^{\alpha_1}\cdot ...\cdot G_m^{\alpha_m} \cdot x_1^{\beta_1}\cdot ... \cdot x_n^{\beta_n} +p\cdot G_1, $$
for some $G_1\in \Z_q[x_1,...,x_n].$

Clearly, $G_1(\Phi^n)=0,$ so $G_1\in Z_q[x_1,...x_n].$   Repeating the above procedure, we get $G_1$ as a linear combination of  $g_1^{\alpha}x^{\beta}$ plus $pG_2,$ and so on. Combining all these together and using the completeness of $\Z_q$ in $p-$adic topology, we prove the lemma.
\endproof

From \cite[Theorem 3.2]{BS}, for some  $q=p^{\tau}$ there is a point
$x=(x_1,...,x_n)\in V(\F_q^n)$ such that $\Phi (x) =x^Q$ for some
$Q=p^l.$ Additionally, we can choose $x$ to lie outside of any fixed
Zariski closed subset. So we choose $x$ to be a smooth point of $V$,
where the restriction of the tangent map $\Phi_*$ to the Zariski
tangent space $T_{x}V$ is injective. This implies that
$\Phi_*|_{T_{x^Q}(V)},$ $\Phi_*|_{T_{x^{Q^2}}(V)},...$ are all
injective. Suppose $N \in \N$ is such that $\Phi ^N(x)=x.$ Then
$\Phi^N_*|_{T_x(V)}$ is injective as a composition of injective
linear operators. So it is invertible.

By our choice of notation, $T_x(V)$ is a vector space over $\F_p^{alg}.$ But it has an $\F_q-$basis, because $x$ is defined over $\F_q$ and $V$ is defined over $\F_p.$ In this basis the matrix of $\Phi^N$ has coefficients in the finite field $\F_q$. It is invertible, so some power of it, $\Phi^M_*$ is identity.

By Theorem 2.4 there exists a point $X\in V(\Z_q)$ such that $\pi (X)=x.$  By definition of $M,$ $\Phi^M(X) \equiv X \mod p .$

The following general observation is very important.
\begin{lemma} Suppose $P=P(x_1,x_2,...,x_n)\in {\Z}_q[x_1,...,x_n].$ Suppose $A=(A_1,...,A_n)\in {\Z}_q^n.$ Suppose $l \in \N.$ Denote by $\nabla P(A)$ the formal gradient of $P,$ evaluated at $A.$ Then for every $Y=(Y_1,...Y_n)\in {\Z}_q^n$ we have the following.
$$P(A+Y\cdot p^{l}) \equiv P(A) + (\nabla P(A) \cdot Y)p^{l} \mod p^{l+1}$$
\end{lemma}

\proof We rewrite the polynomial $P$ as a linear combination of products of powers of $(x_i-A_i).$ When evaluated at $A+Y\cdot p^{l},$ the terms of degree at least two are zero modulo $p^{2l},$ so are zero modulo $p^{l+1}.$  The linear term is exactly the dot product of the gradient and $Y\cdot p^l.$
\endproof

An immediate corollary is the following.

\begin{lemma} Suppose $\Omega$ is an $n-$variable  polynomial automorphism, with coefficients from $\Z_q,$ and $X\in \Z_q^n$ is a point such that $\omega (x)=x,$ where $x=\pi(X)$ and $\omega$ is a reduction of $\Omega$ modulo $p.$ Denote
the iduced map on the tangent space at $x$ by $\omega_*.$ Then for any $Y \in Z_q^n$ with $y=\pi (Y)$ we have
$$\Omega (X+p^lY) = \Omega +p^l\cdot \omega_*(y) \mod p^{l+1} $$
(Note that $\omega_*(y)$ is only defined modulo $p$, but its product with $p^l$ makes sense modulo $p^{l+1}$).
\end{lemma}

\proof
Apply the previous Lemma to each component of $\Omega ,$ with $A=X$.
\endproof

Now we go back to our map $\Phi$ that fixes the subvariety $V$.
By the choice of $M$ and $X,$ $\Phi^M(X)\equiv X \mod p.$ Also, $\phi^M_*$ restricted to the tangent space of $V(\F_q)$ at $x$  is the identity. Denote by $\alpha^{(1)}\in \F_q^n$ the divided difference $\frac{\Phi^M(X)-X}{p} (\mod p)$.

\begin{lemma} For any $X'\equiv X \mod p$ with $X'\in V(\Z_q)$,
$$\Phi ^M (X')  \equiv X' + \alpha^{(1)}\cdot p \mod p^2$$
\end{lemma}

\proof
First of all, $X'=X+pY $ for some $Y .$ We are going to show that $\pi (Y) \in T_*(V) (x).$ This is equivalent to showing that any polynomial in $\F_q[x_1,...x_n]$ that vanishes on $V$ has its gradient, evaluated at $x,$ vanishing at $\pi(y).$ Suppose $g$ is such a polynomial. By Lemma 2.6 there exists $G\in I(V(\Z_q))$ such that $g=\pi (G).$  By Lemma 2.7,
$$G(X')=G(X+pY) \equiv G(X) +p (\nabla G)(X)\cdot Y \mod p^2$$
Since $G(X')=G(X)=0,$ this implies that $\nabla G(X)\cdot Y \equiv 0 \mod p,$
so $\nabla g (x) \cdot y =0$ in $\F_q$.

Now we apply Lemma 2.8. to the map $\Omega = \Phi ^M.$
$$\Phi^M (X') \equiv \Phi ^M (X) + p\Phi^M_* (X) (Y) \equiv (X+p\alpha^{(1)}) +pY \equiv X' + \alpha^{(1)}\cdot p \mod p^2 $$
\endproof

As an immediate corollary, we have the following.
\begin{lemma}  For any $X'\equiv X mod p$ with $X'\in V(\Z_q)$,
$$\Phi ^{pM} (X')  \equiv X' \mod p^2$$
\end{lemma}

\proof We know that $\Phi^M$ fixes $V$. So, by induction, for all natural $j$
$$\Phi ^{jM} (X')  \equiv X' + \alpha^{(1)}\cdot p\cdot j \mod p^2$$
\endproof

Now we denote by $\alpha^{(2)}\in \F_q^n$ the divided difference $\frac{\Phi^{pM}(X)-X}{p^2} (\mod p)$.
We have the following lemma.

\begin{lemma} For any $X'\equiv X mod p^2$ with $X'\in V(\Z_q)$,
$$\Phi ^{pM} (X')  \equiv X' + \alpha^{(2)}\cdot p^2 \mod p^3$$
\end{lemma}

\proof The proof analogous to that of Lemma 9, using Lemmas 2.7 and 2.8. The details are left to the reader.
\endproof

As a corollary of Lemma 2.11 we get that $\Phi ^{p^2M} (X')  \equiv X' \mod p^3$, and so on.
Putting it all together, we get the following theorem.

\begin{theo} \label{MainAlg} Suppose $\Phi$ is an $n-$variable polynomial map with integer coefficients, and $V$ is the Zariski closure of the image of $\Phi ^n$ (over $\Z$). Suppose $W$ is a proper subscheme of $\ V$. Then for every sufficiently large prime $p$ there exist $q=p^l$ and a point $x\in V({\F}_q) \setminus W({\F}_q) $ such that for every $X\in V(\Z_q)$ with $\pi (X)=x$ we have
$$\Phi^{ap^k}\equiv X (\mod p^k),$$
where $a$ is fixed and $k$ is arbitrary.

In particular, the point $X$ is uniformly recurrent for $\Phi$ in the $p-$adic topology on $V(\Z_q)$.
\end{theo}

\section{Residually finite groups}
\label{sect3}

In this section, we shall prove the following

\begin{theo}\label{th4}  Let $\phi$ be any injective endomorphism of a free group $F_k$. Then for every sufficiently large prime $p$, the HNN extension $\HNN_\phi(F_k)$ has a subgroup of finite index that is residually (finite $p$-group) and also is an ascending $\HNN$-extension of a free group.
\end{theo}

\proof Let $F_k=\la x_1,...,x_k\ra$. The endomorphism $\phi$ is
defined by $k$ words $$w_1(x_1,...,x_k), ..., w_k(x_1,...,x_k),$$
the images of $x_1,...,x_k$. For every finite group $G$ consider the
map $\phi_G\colon G^k\to G^k$ defined as follows: $$(g_1,...,g_k)\to
(w_1(g_1,...,g_k), ...,w_k(g_1,...,g_k)).$$ Let $H=\HNN_\phi(F_k)$.

Suppose that a point $(g_1,...,g_k)\in G^k$ is periodic with respect
to the map $\phi_G$ and $l$ is the length of the period. Consider
the wreath product $P$ of $G$ and the cyclic group $C_l=\la c\ra$ of
order $l$, i.e. $P$ is the semidirect product of $G^l$ and $C_l$
where the elements of $C_l$ cyclically permute the factors of $G^l$.
Let $t$ be the free letter of $\HNN_\phi(F_k)$. It was observed in
\cite[page 346]{BS} that then the map \begin{equation}\label{eq01}
t\to c,\quad x_i\to (g_i, \phi_G(g_i),
\phi^2_G(g_i),...,\phi^{l-1}_G(g_i))\end{equation} extends to a
homomorphism from $H$ into $P$.

Consider the ring of $2\times 2$-matrices $M=M(2,\bar\F_p)$ over the
algebraic closure of $\F_p$ as a copy of the affine space of
dimension 4. Replacing the inverses $x^{-1}$ in the words $w_i$ by
symbols $\adj(x)$ interpreted as the adjoint matrix, we turn
$\phi_M$ into a polynomial map $\bar\phi\colon M^k\to M^k$. Assume
that $p$ is large enough (as in Theorem 2.12). By \cite[Theorem
3.2]{BS}, this map has a smooth periodic point $u=(u_1,...,u_k)$ in
the variety $\overline{\phi^{4k}(M^k)}$ (here $\bar{}$ means Zariski
closure) which does not belong to the proper subvariety
$\{(y_1,...,y_k)\in M^k\mid \det(y_1y_2...y_k)=0\}$, that is all the
matrices $u_i$ are non-singular.

Let $\F_q$ be the (finite) field over which all the matrices $p_i$ are defined. We can view $u_i$ as matrices in the projective linear group $U=\PGL(2,\F_q)$. Note that the map induced by $\bar\phi$ on $U$ coincides with $\phi_U$. Hence $(u_1,...,u_k)$ is a smooth periodic point of $\phi_U$.  Hence the variety $V=\overline{\phi_S^{4k}(S^k)}$ contains a smooth periodic point for $\phi_S$.

Consider an unramified at $p$ finite extension $K$ of $\Q$ with the ring of integers $\OO$ and a maximal ideal $(p)$ of $\OO$ with the quotient field $\F_q$. Let, as before, $\Z_q$ be the $p$-adic completion of $\OO$. Then the group $\SL(2,\OO)$ naturally embeds into $\SL(2,\Z_q)$. By \cite[Theorem 4.3]{BG}, there exists a free non-Abelian subgroup $\Gamma$ in $\SL(2,\OO)$ that is dense in the pro-finite topology induced by the congruence subgroups of $\SL(2,\OO)$ modulo $p^k$, $k\ge 0$. Hence the homomorphism $\mu\colon \SL(2,\OO)\to \SL(2,\OO/(p))$ is surjective on $\Gamma$.
Again by \cite[Theorem 4.3]{BG} (see also \cite[Corollary 4.4]{BG}), there exist preimages $\bar u_1,...,\bar u_k$ of $u_1,...,u_k$ that freely generate a free subgroup $F$ of $\SL(2,\OO)$ which we shall identify with $F_k$. Consider the $p$-adic variety $V(\Z_q)$. It contains the subset $\phi_F^{m}(F)$ for all $m\ge 4k$. Since the point $(u_1,...,u_k)$ is periodic, the point $u=\phi_H^m(\bar u_1,...,\bar u_k)$ for a divisible enough $m$ also is a preimage of $(u_1,...,u_k)$ under $\mu$. The coordinates of $u$ also freely generate a free subgroup since $\phi$ is injective. Hence we can assume that $(\bar u_1,...,\bar u_k)$ belongs to $V(\Z_q)$.

Consider the sequence of congruence subgroups of $F$ corresponding to powers of $p$:

$$F>F^{(1)}>F^{(2)}>...>... $$
The intersection $\cap F^{(i)}$ is $\{1\}$ and factor-groups $F^{(i)}/F^{(i+1)}$ are $p$-groups for every $i\ge 1$.
Let $\gamma_i$ be the natural homomorphisms of $F$ onto $\Gamma_i=F/F_i$. By Theorem \ref{MainAlg}, for every $i\ge 1$, the point $(\gamma_i(\bar u_1),...,\gamma_i(\bar u_k))$ in $\Gamma_i^k$ is periodic with respect to $\phi_{\Gamma_i}$ with period $l_i=ap^{i-1}$ for some integer constant $a$ (the length of the period of the point $(u_1,...,u_k)$ from $\Gamma_1^k$).

For every $i$ let $C_{l_i}$ be the cyclic group of order $l_i$. Consider the homomorphism $\nu_i$ from $H$ into the wreath product $P_i=\Gamma_i \wr C_{l_i}$ defined as in (\ref{eq01}). These homomorphisms separate all elements of $\HNN_\phi(F_k)$. Indeed, every non-identity element in $\HNN_\phi(F_k)$ can be represented as $t^nwt^{n'}$ for some $n,n'\in \Z$ and some word $w\in F_k$ where either $w\ne 1$ or $n+n'\ne 0$. If $w\not\in H_i$ for some $i$, then this element is not in the kernel of $\nu_i$. If $w=1$ and $n+n'\ne 0$, then it is not in the kernel of any $\nu_i$ with $i>|n+n'|$.

For every $i\ge 1$, consider the natural homomorphism $\mu_i$ from $\Gamma_i^{l_i}$ to $\Gamma_1^{l_i}$. Note that $\mu_i\nu_i(F)$ is a $k$-generated subgroup of the direct power of $\Gamma_1$, and that the kernel $K_i$  of $\nu_i$ is a $p$-group.  Hence $\mu_i\nu_i(F)$ belongs to the variety of groups generated by the finite group $\Gamma_1$. Since a variety generated by a finite group is locally finite \cite{Neu}, there exists a constant $M$ such that $|\mu_i\nu_i(F)|\le M$ for every $i$. The subgroup $K_i$ is normal in $P_i$, so $K_i\cap \nu_i(F)$ is normal in $\nu_i(H)$. The group $\nu_i(H)/(K_i \cap \nu_i(F)$ is an extension of a group $E_i$ of order at most $M$ by a cyclic group of order  $l_i=ap^{i-1}$ for some constant $a$. The centralizer of $E$ in $\nu_i(H)/(K_i \cap \nu_i(F))$ has index at most $M^M$, hence there exists a constant $M_1$ such that the group $\mu_i\nu_i(H)$ contains a $p$-subgroup of index at most $M_1$. Since the kernel $K_i$ is a $p$-group, the group $\mu_i(H)$ has a $p$-subgroup of index $\le M_1$. Hence $H$ has a subgroup $N$ of index at most $M_1$ that is residually (finite $p$)-group. We can assume that $N$ is a normal subgroup. 

It remains to note that $N$ is generated by the intersection $F_k\cap N$ and $t_1=t^{ap^s}$ for some $s$. Moreover, $t_1(F_k\cap N)t_1\iv\subseteq F_k\cap N$ since $t_1F_kt_1\iv\subseteq F_k$ and $t_1\in N$. Hence $N$ is an ascending $\HNN$-extension of a free group (see, for instance, \cite{DS}).
\endproof

\end{document}